\documentclass[12pt]{article}
\usepackage{amsmath,amsthm,amssymb,amscd}
\usepackage{latexsym}
\newtheorem{thm}{Theorem}[section]
 
 \newtheorem{lem}[thm]{Lemma}
 \newtheorem{prop}[thm]{Proposition}
 \theoremstyle{definition}
 
 \theoremstyle{remark}

 \numberwithin{equation}{section}

\title
{K\"{a}hler-Ricci flow on  homogeneous toric bundles}

\author{ Hong Huang}
\date{}
\begin{document}
\maketitle
\begin{abstract}
  Assume that $X$ is a homogeneous toric bundle of the form $G^{\mathbb{C}}\times_{P,\tau} F$ and is Fano, where $G$ is a compact semisimple Lie group with  complexification $G^\mathbb{C}$,  $P$ a parabolic subgroup  of $G^\mathbb{C}$, $\tau:P\rightarrow (T^m)^\mathbb{C}$ is a surjective homomorphism from $P$ to the algebraic torus $(T^m)^\mathbb{C}$, and $F$ is a compact toric manifold of complex dimension $m$. In this note we show that the normalized K\"{a}hler-Ricci flow on $X$ with a  $G\times T^m$-invariant initial K\"{a}hler form in $c_1(X)$ converges, modulo the algebraic torus action, to a K\"{a}hler-Ricci soliton. This extends a previous work of X. H. Zhu. As a consequence we recover a result of Podest\`{a}-Spiro.



\end{abstract}
\maketitle


\section {Introduction}

In a seminal work \cite{WZ} Wang and Zhu proved the existence of a K\"{a}hler-Ricci soliton on any toric Fano manifold.
This result was recovered and generalized in some later works. Here we only mention two such works of particular relevance to this note: Podest\`{a} and Spiro \cite{PS} generalized Wang-Zhu's result to the case of homogeneous toric bundles  which are Fano, and Zhu \cite{Z} recovered his result with Wang using  K\"{a}hler-Ricci flow. In this  note we extend \cite{Z} to the case of  Fano homogeneous toric bundles and recover Podest\`{a} and Spiro's result in \cite{PS}  via K\"{a}hler-Ricci flow. More precisely we show

\begin{thm} \label{thm 1.1} \ \
  Assume that $X$ is a homogeneous toric bundle of the form $G^{\mathbb{C}}\times_{P,\tau} F$ and is Fano,  where $G$ is a compact semisimple Lie group with  complexification $G^\mathbb{C}$,  $P$ a parabolic subgroup  of $G^\mathbb{C}$, $\tau:P\rightarrow (T^m)^\mathbb{C}$ is a surjective homomorphism from $P$ to the algebraic torus $(T^m)^\mathbb{C}$, and $F$ is a compact toric manifold of complex dimension $m$. The normalized K\"{a}hler-Ricci flow on $X$ with a  $G\times T^m$-invariant initial K\"{a}hler form in $c_1(X)$ converges, modulo the algebraic torus action, to a K\"{a}hler-Ricci soliton.
\end{thm}

We emphasize that in Theorem 1.1 we do not assume the existence of a K\"{a}hler-Ricci soliton on the Fano homogeneous toric bundle $X$; the existence of K\"{a}hler-Ricci soliton on $X$ is part of the conclusion (which, however, was previously known  by \cite{PS}, as mentioned above).  On the other hand,  using a result in \cite{TZZZ} Dervan and Sz\'{e}kelyhidi \cite{DS} prove that if a  Fano manifold admits a K\"{a}hler-Ricci soliton $\omega_{KRS}$, then the normalized K\"{a}hler-Ricci flow with any initial K\"{a}hler form in the first Chern class, converges to $\omega_{KRS}$ modulo the  action of automorphisms.  So  the following is true: Assume that $X$ is a homogeneous toric bundle of the form $G^{\mathbb{C}}\times_{P,\tau} F$ and is Fano, then
 the normalized K\"{a}hler-Ricci flow on $X$ with any initial K\"{a}hler form in $c_1(X)$ converges, modulo the  action of automorphisms, to a K\"{a}hler-Ricci soliton.    I do not know how to prove this fact directly. However, note that this fact does not imply the whole  Theorem 1.1: the conclusion of Theorem 1.1 is slightly stronger; of course, the condition of Theorem 1.1 is also slightly stronger.

There are many beautiful works on the K\"{a}hler-Ricci flow, see for examples \cite{BEG},  \cite{T}, \cite{To} and \cite{W}  for some recent surveys.

 The proof of Theorem 1.1  follows closely the lines of \cite{Z}, and we also borrow  some key results from \cite{PS} (so the proof in this note is not  completely independent of that in \cite{PS}). In \cite{WZ}, to obtain some of the key estimates Wang and Zhu converted the K\"{a}hler-Ricci soliton equation - a complex Monge-Amp\`{e}re equation - on  a toric Fano manifold to a real Monge-Amp\`{e}re equation on the Euclidean space. Similarly, one of the key steps in \cite{PS} (see Proposition 5.2 in \cite{PS}) is to convert the K\"{a}hler-Ricci soliton equation on a  homogeneous toric bundle which is Fano to an equation on the open dense  orbit of the $(T^m)^\mathbb{C}$-action on the fiber, which is actually a real Monge-Amp\`{e}re equation on $\mathbb{R}^m$.

 In the flow case, first Zhu  \cite{Z} converted the K\"{a}hler-Ricci flow - a parabolic complex Monge-Amp\`{e}re equation (PCMAE)- on a toric Fano  manifold to a parabolic real Monge-Amp\`{e}re equation on the Euclidean space. Then he adapted some of the key estimates on the real Monge-Amp\`{e}re equation in \cite{WZ} to the parabolic case with the help of a deep estimate of Perelman (see \cite{ST}).

 Here (see Section 2.3) using some results in \cite{PS} we also convert the  K\"{a}hler-Ricci flow on a  homogeneous toric bundle which is Fano to a
parabolic real Monge-Amp\`{e}re equation on the Euclidean space, which is very similar to the one in \cite{Z}; the only difference is that in our case there is an extra  term in the equation, which turns out to be  harmless due to a result in \cite{PS}.  Note that in Section 3.4 of \cite{Do} Donaldson made some  modifications  to the original approach in \cite{WZ}. In the first step of   the estimates for the parabolic real Monge-Amp\`{e}re equation we follow  Donaldson's  modifications (see Section 3).  We also observe that under the normalized K\"{a}hler-Ricci flow on a Fano homogeneous toric bundle with a  $G\times T^m$-invariant initial K\"{a}hler form in $c_1(X)$, the volume of any fiber is a constant.

In Section 2.1, we briefly recall some basic concepts about toric Fano manifolds, and give  alternative proofs of some known facts related to the K\"{a}hler potentials and moment maps of  toric Fano manifolds (see the proofs of Facts  4 and 5 there).  In Section 2.2, we give a brief introduction  to homogeneous toric bundles following \cite{PS}.  In Section 2.3 we reduce the normalized K\"{a}hler-Ricci flow on a homogeneous toric bundle which is Fano to a parabolic real Monge-Amp\`{e}re equation on the Euclidean space. In Section 3 we prove our Theorem 1.1 following \cite{Z} in most steps, except in the first step where we mainly follow (and expand) part of Section 3.4 in \cite{Do} instead.

\vspace*{0.4cm}

\section{Reduction of K\"{a}hler-Ricci flow to PMAE}

\vspace*{0.4cm}

First we fix some conventions and notations. For a compact K\"{a}hler manifold $X$ of complex dimension $n$  with K\"{a}hler metric $g=\sum_{i,j=1}^n g_{i\bar{j}}dz^i\otimes d\bar{z}^j$,  the corresponding K\"{a}hler form $\omega=\frac{\sqrt{-1}}{2\pi}\sum_{i,j=1}^n g_{i\bar{j}}dz^i\wedge d\bar{z}^j$.   Let $\text{Ric}(\omega)=\frac{\sqrt{-1}}{2\pi}\sum_{i,j=1}^n R_{i\bar{j}}dz^i\wedge d\bar{z}^j$
be the Ricci form of $\omega$, where
\begin{equation*}
R_{i\bar{j}}=-\frac{\partial^2 }{\partial z^i \partial \bar{z}^j}\log \det (g_{k\bar{l}}).
\end{equation*}
$\text{Ric}(\omega)$ represents the first Chern class $c_1(X)$.

Let $G$ be a Lie group with Lie algebra $\mathfrak{g}$. Suppose  $G$ acts on a  manifold $X$. For any $Y \in \mathfrak{g}$,  we will denote the induced vector field on $X$ by $\hat{Y}$ (see line 3 on p. 122 of \cite{H}). We also denote the set of $G$-principal points in $X$ by $X_{\text{reg}}$.

\vspace*{0.4cm}

\noindent {\bf 2.1  Toric Fano manifolds}

\vspace*{0.4cm}

We give a brief review of toric manifolds following \cite{BS} and \cite{M}. One can also consult \cite{CLS}, \cite{F} and \cite{O}  for more details. Let $M$ be a free abelian group of rank $m$, $N=$ Hom$(M,\mathbb{Z})$, $M_\mathbb{R}:=M\otimes_\mathbb{Z}\mathbb{R}$, $N_\mathbb{R}:=N\otimes_\mathbb{Z}\mathbb{R}$. Denote by $\langle\cdot, \cdot\rangle: M_\mathbb{R} \times N_\mathbb{R}\rightarrow \mathbb{R}$ the natural pairing.
 Let $F$ be a smooth projective toric $m$-fold defined by a complete fan $\Sigma$ of regular cones $\sigma \subset N_\mathbb{R}$.  We have a maximal (algebraic) torus $(T^m)^\mathbb{C}=\text{Hom} (M, \mathbb{C}^*) \subset$ Aut$(F)$, $(T^m)^\mathbb{C} \cong (\mathbb{C}^*)^m$. $(T^m)^\mathbb{C}$ has an open dense orbit $U\subset F$, $U\cong (T^m)^\mathbb{C}$. Using the above notation $U=F_{\text{reg}}$. We choose a basis for $M$ and the dual basis for $N$, so we can identify $ M_\mathbb{R}$ ($N_\mathbb{R}$) with $\mathbb{R}^m$ and $(T^m)^\mathbb{C}$ with $(\mathbb{C}^*)^m$. We also fix a point in $U$ so we can also identify $U$ with $(T^m)^\mathbb{C}$ (and also with $(\mathbb{C}^*)^m$). For $z=(z_1,\cdot\cdot\cdot, z_m) \in (\mathbb{C}^*)^m$  let $t_i=\log |z_i|^2$.
To $a=(a_1,\cdot\cdot\cdot,a_m)\in M$ we associate the algebraic character $\chi^a:(T^m)^\mathbb{C}\rightarrow \mathbb{C}^*$,
\begin{equation*}
\chi^a(z):=z_1^{a_1}\cdot\cdot\cdot z_m^{a_m}.
\end{equation*}
Then
\begin{equation*}
|\chi^a(z)|^2=e^{\langle a,\underline{t}\rangle},
\end{equation*}
where $\underline{t}=(t_1,\cdot\cdot\cdot,t_m)\in N_\mathbb{R}$ and  $t_i=\log |z_i|^2$ as above.

Now we assume that the toric manifold $F$ is Fano. Then by Demazure (see Theorem 2.1 in \cite{M}) we get an $m$-dimensional compact convex polytope $\Delta \subset M_\mathbb{R}$ defined by the inequalities $\langle \cdot, f\rangle \leq 1$, where $f$ runs over  primitive integral generators of all
1-dimensional cones in $\Sigma$. ($\Sigma$ is called the normal fan of $\Delta$.) Let $\{p^{(0)},\cdot\cdot\cdot,p^{(s)}\}=M\cap \Delta$. By Demazure (see  Section 2.3 in \cite{O}) we have the anticanonical embedding $F\hookrightarrow \mathbb{CP}^s$ defined by $\chi^{p^{(0)}},\cdot\cdot\cdot,  \chi^{p^{(s)}}$. Define a function
$v^0: U \rightarrow \mathbb{R}$ by
\begin{equation*}
v^0(z)=\log(\sum_{k=0}^s|\chi^{p^{(k)}}(z)|^2).
\end{equation*}
Following the convention on p. 711 in \cite{M}, we still denote by $v^0$ the function $N_\mathbb{R}\rightarrow \mathbb{R}$ giving by
\begin{equation*}
\underline{t}\mapsto \log(\sum_{k=0}^se^{\langle p^{(k)},\underline{t}\rangle}).
\end{equation*}
So we have $v^0(z)=v^0(\underline{t})$, where $z=(z_1, \cdot\cdot\cdot, z_m)\in U$,   $\underline{t}=(t_1,\cdot\cdot\cdot,t_m)\in N_\mathbb{R}$, and  $t_i=\log |z_i|^2$ as above, and we may view $v^0$ either as a function defined on $U$ or as a function on $N_\mathbb{R}$. Note that  the K\"{a}hler form $\frac{\sqrt{-1}}{2\pi}\partial \bar{\partial}v^0$ on $U$ extends to a K\"{a}hler form on $F$, denoted by $\hat{\omega}_0$, which is the pull-back of the  K\"{a}hler form $\frac{\sqrt{-1}}{2\pi}\partial \bar{\partial} \log \sum_{k=0}^s|z_k|^2$  (corresponding to the Fubini-Study metric) on $\mathbb{CP}^s$ via the anticanonical embedding $F\hookrightarrow \mathbb{CP}^s$ defined by $\chi^{p^{(0)}},\cdot\cdot\cdot,  \chi^{p^{(s)}}$.

Let $z_i=e^{w_i}$, where $w_i=\frac{1}{2}t_i+\sqrt{-1}\theta_i$ (in particular, $t_i=\log |z_i|^2$). Then we have $\frac{\partial}{\partial w_i}=z_i\frac{\partial}{\partial z_i}$, $\frac{\partial}{\partial \bar{w}_i}=\bar{z}_i\frac{\partial}{\partial \bar{z}_i}$, $dw_i=\frac{1}{z_i} dz_i$,  and $d\bar{w}_i=\frac{1}{\bar{z}_i} d\bar{z}_i$; moreover, $\frac{\partial}{\partial w_i}=\frac{1}{2}(\frac{\partial}{\partial \frac{t_i}{2}}-\sqrt{-1}\frac{\partial}{\partial \theta_i})$,
$\frac{\partial}{\partial \bar{w}_i}=\frac{1}{2}(\frac{\partial}{\partial \frac{t_i}{2}}+\sqrt{-1}\frac{\partial}{\partial \theta_i})$, $dw_i=\frac{1}{2}dt_i+\sqrt{-1}d\theta_i$, and $d\bar{w}_i=\frac{1}{2}dt_i-\sqrt{-1}d\theta_i$.
As $v^0: U \rightarrow \mathbb{R}$ is $T^m$-invariant, we  have
 \begin{equation*}
  \partial\bar{\partial}v^0
=  \frac{\partial^2 v^0}{\partial w_i \partial \bar{w}_j} dw_i\wedge d\bar{w}_j
= \frac{\partial^2 v^0}{\partial t_i \partial t_j}dw_i\wedge d\bar{w}_j
=-\sqrt{-1}\frac{\partial^2 v^0}{\partial t_i \partial t_j}dt_i\wedge d\theta_j.
 \end{equation*}
  Sometimes we'll denote $\frac{\partial^2 v^0}{\partial t_i \partial t_j}$ by $v_{ij}^0$.

Given a toric Fano manifold $F$ as above, for any $T^m$-invariant K\"{a}hler form $\rho \in c_1(F)$, by Calabi-Yau theorem \cite{Yau} (one can find  a much simpler proof in our toric case in \cite{Do}) there is a unique K\"{a}hler form $\omega_\rho$ with $\text{Ric}(\omega_\rho)=\rho$. Consider the unique moment map
$\mu_\rho:F\rightarrow \mathfrak{t}^*$ relative to $\rho$ (that is, $d\langle \mu_\rho(\cdot), \xi\rangle=-\iota_{\hat{\xi}}\rho$ for any  $\xi \in \mathfrak{t}$) with $\int_F\mu_\rho \omega_\rho^m=0$, which is called metrically normalized (following \cite{PS}; compare Remark 9.4.2 in \cite{M}). By Atiyah and Guillemin-Sternberg $\mu_\rho(F)$  is a compact convex polytope.  It is pointed out in Section 2.3 in \cite{PS} that $\mu_\rho(F)$ does not depend on the choice of $\rho$, which also follows from Facts 2 and 5 below.   Denote $\mu_\rho(F)$ by $\Delta_F$.

Note that the set of all vertices of $\Delta$ is contained in $M \cap \Delta$.  By Demazure (see Theorem 2.1 in \cite{M}), the set of all  vertices of $\Delta$,  denoted by $\{p^{(1)}, \cdot\cdot\cdot, p^{(l)}\}$,  corresponds  one-to-one  to the set of all  $m$-dimensional cones in the fan $\Sigma$,  denoted by $\{\sigma_1, \cdot\cdot\cdot, \sigma_l\}$,  with $\langle p^{(k)}, f\rangle=1$ for all fundamental generators $f$ of $\sigma_k$, $k=1, 2, \cdot\cdot\cdot, l$.
Note that $N_\mathbb{R}=\cup_{k=1}^l \sigma_k$. Let $\bar{v}(\underline{t}):=\max_{1\leq k \leq l} \langle p^{(k)},\underline{t}\rangle$.  It is easy to see that $\bar{v}$ is a piecewise linear function with $\bar{v}(\underline{t})=\langle p^{(k)}, \underline{t}\rangle$ for $\underline{t}\in \sigma_k$.

The following five facts are known and will be used later.

\vspace*{0.4cm}

\noindent {\bf Fact 1} (cf. \cite{BS}).   We have $0 \leq v^0-\bar{v} \leq \log (s+1)$.

\vspace*{0.4cm}

\noindent {\bf Proof}.  Compare \cite{BS}. Clearly $\bar{v}(\underline{t})=\max_{0\leq k \leq s} \langle p^{(k)},\underline{t}\rangle$.  Now
\begin{equation*}
e^{\bar{v}(\underline{t})} \leq \sum_{k=0}^se^{\langle p^{(k)},\underline{t}\rangle} \leq (s+1) e^{\bar{v}(\underline{t})},
\end{equation*}
and the desired result follows.    \hfill{$\Box$}

\vspace*{0.4cm}

Let $\omega \in c_1(F)$ be a $T^m$-invariant K\"{a}hler form.  By Section 3 of \cite{M} (compare Theorem 4.3 in \cite{G}) there exists a $T^m$-invariant function
$u \in C^\infty(F_{\text{reg}})$  such that $(\sqrt{-1})^m e^{-u}\prod_{i=1}^m(dw_i\wedge d\bar{w}_i)$ extends to a volume form on $F$ and $\omega|_{F_{\text{reg}}}=\frac{\sqrt{-1}}{2\pi}\partial \bar{\partial}u$. We also view $u$ as a function on $N_\mathbb{R}\cong \mathbb{R}^m$ as we do with $v^0$.

\vspace*{0.4cm}

\noindent {\bf Fact 2}  (Atiyah and Guillemin-Sternberg, cf. Theorem 4.2 in \cite{M}).  Let $u$ be as above. The map $\phi(z):=Du( \underline{t})$ is the restriction to $U$ of a moment map relative to  $2\pi\omega$ for the $T^m$-action on $F$, where $z=(z_1,\cdot\cdot\cdot,z_m)\in U$, $\underline{t}=(t_1, \cdot\cdot\cdot,t_m)\in N_\mathbb{R}$ with $t_i=\log |z_i|^2$. Moreover
the gradient map
\begin{equation*}
Du: N_\mathbb{R}\rightarrow M_\mathbb{R}(\cong \text{Lie}(T^m)^*),        \hspace{4mm}   \underline{t}\mapsto (\frac{\partial u}{\partial t_1}(\underline{t}),\cdot\cdot\cdot, \frac{\partial u}{\partial t_m}(\underline{t}))
\end{equation*}
is a diffeomorphism from $N_\mathbb{R}$ to $\Delta \setminus \partial \Delta$.

\vspace*{0.4cm}

\noindent {\bf Proof}.  In fact  Theorem 4.2 in \cite{M} is slightly stronger; compare Remark 4.3 there.   For more details see the proof of Theorem 8.2 in \cite{M}
(compare also  for example,  Exercise 12.2.8 in \cite{CLS} for the case $u=v^0$).
Let's check the first statement in Fact 2.
 Suppose $\frac{\partial}{\partial \theta_i}$ is induced by  $Y_i \in \mathfrak{t}=\text{Lie}(T^m)$. Then $\phi(z)=\sum_{i=1}^m\phi_i(z)Y_i^*$, where $\phi_i(z)=\frac{\partial u}{\partial t_i}$. We have
\begin{equation*}
\begin{split}
  \iota_{\frac{\partial}{\partial \theta_i}} \sqrt{-1} \partial  \bar{\partial}u
=& \iota_{\frac{\partial}{\partial \theta_i}}   \frac{\partial^2 u}{\partial t_j \partial t_k}dt_j\wedge d\theta_k\\
=& -\frac{\partial^2 u}{\partial t_i \partial t_j}dt_j \\
=& -d(\frac{\partial u}{\partial t_i}) \\
=& -d\phi_i,
\end{split}
\end{equation*}
 and the first statement in Fact 2 follows.

For the second statement in Fact 2 we only briefly recall some of the key steps in the proof of Theorem 8.2 in \cite{M}  and give an alternative argument under an extra condition  for one of them.
On p. 720-721 of \cite{M} it is shown that $\phi$ naturally extends to $F$, which will be denoted by $\tilde{\phi}$. On p. 722 of \cite{M} it is shown that all vertices of $\Delta$ are in the image of $\tilde{\phi}$; below we'll give a simple proof of this result under the extra condition that $u=v^0+\varphi$, where  $\varphi$ is a $T^m$-invariant smooth function  on $F$.

Let $u=v^0+\varphi$ as above, and $u_\lambda(\underline{t})=\lambda^{-1}u(\lambda \underline{t})$ for $\lambda >0$. By Fact 1, $|u_\lambda-\bar{v}|\leq C\lambda^{-1}$, where the constant $C$ depends on the $C^0$-bound of $\varphi$ on $F$.
 Let $\lambda=1, 2, \cdot\cdot\cdot$, and consider the sequence $u_1, u_2, \cdot\cdot\cdot$. Then we have $\text{lim}_{i \rightarrow \infty}u_i=\bar{v}$. (Compare p. 39 in \cite{Do}.) By Theorem 25.7 in \cite{R},
 \begin{equation*}
 \text{lim}_{i \rightarrow \infty} Du(i\underline{t})=\text{lim}_{i \rightarrow \infty} Du_i(\underline{t})=D\bar{v}(\underline{t})=p^{(k)},   \hspace{4mm}      \forall \underline{t} \in \text{int} (\sigma_k),  \hspace{4mm}  k=1,2, \cdot\cdot\cdot,l,
\end{equation*}
where $\text{int} (\sigma_k)$ is the interior of $\sigma_k$.  It follows that $p^{(k)} \in   \overline{Du(\text{int} (\sigma_k))}   \subset    \overline{Du(N_\mathbb{R})}$.

 On p. 723 of \cite{M} it is shown that $\tilde{\phi}(F)\subset \Delta$; compare the proof of Theorem 3.4 in  \cite{BS}, where it is also shown that $D(v^0+\varphi)( N_\mathbb{R}) \subset \Delta$ by using Fact 1.
        \hfill{$\Box$}

\vspace*{0.4cm}

\noindent {\bf Fact 3} (cf. \cite{G}, \cite{A}, p. 38-39 in \cite{Do} and p. 327 in \cite{Z}) Let $u$ be as above. Then the function $u-v^0$  can be extended to a $T^m$-invariant smooth function on $F$, and in particular, $|u-v^0|$ is bounded.

\vspace*{0.4cm}

\noindent {\bf Proof}.  The result is implicitly contained in \cite{G} and \cite{A}. Since both $\omega$  and $\hat{\omega}_0$  are in   $c_1(F)$, there is a $T^m$-invariant smooth function, denoted  by $\varphi$, on $F$ such that
\begin{equation*}
\frac{\sqrt{-1}}{2\pi}\partial \bar{\partial}u-\frac{\sqrt{-1}}{2\pi}\partial \bar{\partial}v^0=\frac{\sqrt{-1}}{2\pi}\partial \bar{\partial}\varphi
\end{equation*}
on $F_{\text{reg}}$.

By Fact 2 we have Im$(Du)=\Delta \setminus \partial \Delta=$Im$(Dv^0)$. Then by Fact 1 and convexity  we have that $\bar{v}-u$ is bounded below (compare the proof of Lemma 3.4 in \cite{WZ}). Combining this with Fact 1 we get that $v^0-u$ is bounded below.  Then we have  that  $v^0-u+\varphi$ is a  harmonic function on $\mathbb{R}^m$ which is bounded below, so it must be a constant.    \hfill{$\Box$}

\vspace*{0.4cm}

It follows that $u$ as above is uniquely determined by $\omega$ up to an additive constant.

\vspace*{0.4cm}

Let $\omega \in c_1(F)$ be a $T^m$-invariant K\"{a}hler form as above. Then there exists a $T^m$-invariant smooth function $\varphi$ on $F$ such that $\omega=\hat{\omega}_0+\frac{\sqrt{-1}}{2\pi}\partial \bar{\partial}\varphi$.
Note that $\varphi$ is uniquely determined up to an additive constant.

\vspace*{0.4cm}

\noindent {\bf Fact 4} (cf. \cite{WZ} and Lemma 4.3 in \cite{S} for the case $\varphi=0$). We have
\begin{equation*}
\sup_{N_\mathbb{R}}|\log\det(\frac{\partial^2 (v^0+\varphi)}{\partial t_i \partial t_j})+v^0+\varphi| < \infty.
\end{equation*}
{\bf Proof}.  As recalled above, the Fano condition implies that we can suppose that  the polytope $\Delta$ is defined by the inequalities $l_r(\cdot)£º=\langle \cdot , \lambda_r\rangle+1\geq 0$,  where $\lambda_r \in N$, $r=1,\cdot\cdot\cdot,d$. (We may view the $\lambda_r$ as inward-pointing normals to codimension 1 faces of $\Delta$.) Let $G_0(x)$ be the Legendre transform of $(v^0+\varphi)(\underline{t})$. By Guillemin \cite{G} we have
\begin{equation*}
(v^0+\varphi)(\underline{t})=\sum_{r=1}^d(-\log l_r(x)+l_r(x))+h_1(x),
\end{equation*}
 and
\begin{equation*}
 G_0(x)=\sum_{r=1}^d l_r(x)\log l_r(x)+h_2(x),
\end{equation*}
where  $x=D(v^0+\varphi)(\underline{t}) \in \Delta \setminus \partial \Delta$ (by Fact 2), and $h_1$ and $h_2$ are smooth functions on (an open neighborhood of) $\Delta$.

Now we have that
\begin{equation*}
\det(\frac{\partial^2 (v^0+\varphi)}{\partial t_i \partial t_j})=(\det \text{Hess} \hspace{1mm} G_0(x))^{-1}=\delta_\Delta(x)\cdot \prod_{r=1}^dl_r(x),
\end{equation*}
where $\delta_\Delta$ is  smooth and positive on $\Delta$. One way to get the last equality is as follows.
Since $\Delta$ satisfies  the Delzant condition \cite{D}, after an affine transformation whose linear part is  given by an element in GL$(m,\mathbb{Z})$ we may assume that the origin  is a vertex of $\Delta$ and near the origin the polytope $\Delta$ is  given by the inequalities $x_1\geq 0, \cdot\cdot\cdot, x_m \geq 0$. Then near the origin we have $G_0(x)=\sum_{i=1}^mx_i\log x_i+v(x)$, where $v$ is smooth.  Similarly we can treat a boundary point which lies on a codimension $i$ face of $\Delta$; see p. 39 of \cite{Do}. Then we can compute with the help of such local expressions of $G_0(x)$ to derive the last equality.   Compare Section 2.4 and Appendix of \cite{A}.  Now  Fact 4 follows.    \hfill{$\Box$}

\vspace*{0.4cm}

\noindent {\bf Fact 5} (cf. \cite{PS}).      The map $\mu_\omega(z):=\frac{1}{2\pi}D(v^0+\varphi)(\underline{t})$ is the restriction to $U$ of the metrically normalized moment map relative to the K\"{a}hler form $\omega$ for the $T^m$-action on $F$, where $z=(z_1,\cdot\cdot\cdot,z_m)\in U$, $\underline{t}=(t_1, \cdot\cdot\cdot,t_m)\in N_\mathbb{R}$ with $t_i=\log |z_i|^2$.

\vspace*{0.4cm}

\noindent {\bf Proof}.  By the first statement in Fact 2 it remains to  verify that $\mu_\omega$ is metrically normalized. One can  show this fact by  inspecting the proof of Lemma 5.1 in \cite{PS}. We give a direct proof below.
   Choose $\psi\in c_1(F)$ such that Ric$(\psi)=\omega$.  Let $h_\omega$ be a Ricci potential of $\omega$, that is, $\text{Ric}(\omega)-\omega=\frac{\sqrt{-1}}{2\pi}\partial \bar{\partial}h_\omega$.   Then we have
\begin{equation*}
-\frac{\sqrt{-1}}{2\pi}\partial \bar{\partial}\log \det ((v^0+\varphi)_{ij})-\frac{\sqrt{-1}}{2\pi}\partial \bar{\partial}(v^0+\varphi) = \frac{\sqrt{-1}}{2\pi}\partial \bar{\partial}h_\omega
\end{equation*}
on $U$. Using Fact 4 and the maximum principle we see that there is a constant $C$ such that
\begin{equation*}
e^{h_\omega}=e^Ce^{-(v^0+\varphi)}\det((v^0+\varphi)_{ij})^{-1}.
\end{equation*}
On the  other hand it is easy to see that there is a constant $C'$ such that
\begin{equation*}
\psi^m=C'e^{h_\omega}\omega^m.
\end{equation*}
(See p. 339-340 of \cite{Yau}.)

Note that on $U$ we have $(2\pi\omega)^m=\det ((v^0+\varphi)_{ij})dt_1\wedge \cdot\cdot\cdot \wedge dt_m \wedge d\theta_1 \wedge \cdot\cdot\cdot \wedge d\theta_m$. Moreover, by Fact 1, $v^0+\varphi \rightarrow +\infty$ as $ \underline{t} \rightarrow \infty$.  Write $\mu_i(z)=\frac{1}{2\pi}\frac{\partial (v^0+\varphi)}{\partial t_i}$. Now
\begin{equation*}
\begin{split}
 \int_U\mu_i\psi^m
=&  C'\int_U\mu_i e^{h_\omega}\omega^m       \\
=&  C''\int_{N_\mathbb{R}}\frac{\partial (v^0+\varphi)}{\partial t_i}e^{-(v^0+\varphi)}dt_1\cdot\cdot\cdot dt_m       \\
=&  -C''\int_{N_\mathbb{R}}\frac{\partial e^{-(v^0+\varphi)}}{\partial t_i}dt_1\cdot\cdot\cdot dt_m     \\
=& 0,
\end{split}
\end{equation*}
and we are done.    \hfill{$\Box$}

\vspace*{0.4cm}

 From Facts 2 and 5 we have $\mu_\omega(F)=\frac{1}{2\pi} \Delta$, the $\frac{1}{2\pi}$-dilation of $\Delta$ w.r.t. the origin.  So $\Delta_F=\frac{1}{2\pi} \Delta$.

\vspace*{0.4cm}

\noindent {\bf 2.2  Homogeneous toric bundles}

\vspace*{0.4cm}

 Now we will follow \cite{PS} to give a brief introduction to homogeneous toric bundles, and refer to \cite{PS} for more details; one can also consult \cite{H} for basics of Lie groups, Lie algebras and symmetric spaces.

 Assume  $G$ is a compact semisimple Lie group. Let $\mathfrak{B}$ be the Cartan-Killing form on $\mathfrak{g}$, we may  view $-\mathfrak{B}$ as an inner product on $\mathfrak{g}$. Let $\mathfrak{g}^\mathbb{C}$ be the complexification of $\mathfrak{g}$. Given a root system $R$ w.r.t. a fixed maximal torus, choose the root vectors  $E_\alpha \in \mathfrak{g}^\mathbb{C}$  corresponding to  roots $\alpha \in R$ via Chevalley's normalization; in particular,  $[E_\alpha, E_{-\alpha}]=H_\alpha$, where $H_\alpha$ is the $\mathfrak{B}$-dual of $\alpha$.

Now let $G$ be a compact semisimple Lie group, and $(F^m, J_F)$ be a compact toric manifold with a complex structure $J_F$ and a holomorphic action of $T^m(\cong (S^1)^m)$. As in \cite{PS} we consider a homogeneous toric bundle $X$ over a generalized flag manifold $V=G^\mathbb{C}/P=G/K$ with fiber $F$,
 \begin{equation*}
X:=G^{\mathbb{C}}\times_{P,\tau} F=G\times_{K,\tau}F,
\end{equation*}
where  $P$ is a parabolic subgroup  of $G^\mathbb{C}$, $\tau:P\rightarrow (T^m)^\mathbb{C}$ is a surjective homomorphism from $P$ to the algebraic torus $(T^m)^\mathbb{C}(\cong (\mathbb{C}^*)^m)$, and $K=P\cap G$.  Note that $V$ is a complex homogeneous space with a natural $G^\mathbb{C}$-invariant complex structure (which will be denoted by $J_V$).  Let $\pi: X\rightarrow V$ be the bundle projection. As in \cite{PS} we identify $F$ with the fiber $F_{eK}:=\pi^{-1}(eK)$.
There is a unique $G^\mathbb{C}$-invariant complex structure (denoted by $J$) on $X$ such that the map $\pi$ is holomorphic and the restriction of $J$ to the fiber $F_{eK}$ is $J_F$.

Note that  we have a holomorphic action of $G^\mathbb{C}\times (T^m)^\mathbb{C}$ on $(X,J)$,
\begin{equation*}
(g',h)\cdot [g,x]:=[g'g, h\cdot x],   \hspace{4mm}      \forall     g', g  \in G^\mathbb{C} \hspace{2mm},  \hspace{2mm}  h \in (T^m)^\mathbb{C}   \hspace{2mm} \text{and} \hspace{2mm} x \in F,
\end{equation*}
compare \cite{Y}.

Recall that $\mathfrak{g}$ has an Ad($K$)-invariant decomposition $\mathfrak{g}=\mathfrak{k}\oplus \mathfrak{m}$. For any fixed Cartan subalgebra (CSA) $\mathfrak{h}\subset \mathfrak{\mathfrak{k}}^\mathbb{C}$ of $\mathfrak{g}^\mathbb{C}$, the associated root system $R$ has a corresponding decomposition $R=R_o+R_\mathfrak{m}$ with $E_\alpha \in \mathfrak{k}^\mathbb{C}$ for $\alpha \in R_o$ and $E_\alpha \in \mathfrak{m}^\mathbb{C}$ for $\alpha \in R_\mathfrak{m}$. The natural $G^\mathbb{C}$-invariant complex structure $J_V$ of $V=G^\mathbb{C}/P$ induces a splitting $R_\mathfrak{m}=R_\mathfrak{m}^+\cup R_\mathfrak{m}^-$, such that
$\sum_{\alpha \in R_\mathfrak{m}^+} \mathbb{C}E_\alpha$ (resp. $\sum_{\alpha \in R_\mathfrak{m}^-} \mathbb{C}E_\alpha$) is the $J_V$-holomorphic (resp. $J_V$-antiholomorphic) subspace $\mathfrak{m}^{(1,0)}$ (resp. $\mathfrak{m}^{(0,1)}$) of $\mathfrak{m}^\mathbb{C}$.
Let
\begin{equation*}
Z_V=-\frac{1}{2\pi}\sum_{\alpha\in R_\mathfrak{m}^+} \sqrt{-1}H_\alpha.
\end{equation*}
Let $Z(K)$ be the center of $K$ and $\mathfrak{z}(\mathfrak{k})$ the Lie algebra of $Z(K)$.  Let $\mathfrak{t}_G:=(\text{ker} \hspace{1mm} d\tau_e)^\bot \cap \mathfrak{z}(\mathfrak{k})$. (We can identify $\mathfrak{t}_G$ with $\mathfrak{t}$ via $d\tau_e$.)
Choose a $-\mathfrak{B}$-orthonormal basis $(Z_1,\cdot\cdot\cdot,Z_m)$ of $\mathfrak{t}_G$ so that $\text{exp}(\mathbb{R}\cdot Z_j)$ is closed for each $j=1,\cdot\cdot\cdot,m$.

 In \cite{PS} Podest\`{a} and Spiro give a criterion for a homogeneous toric bundle $X$ to be Fano, see Theorem A in \cite{PS}, which says that a homogeneous toric bundle $X=G^{\mathbb{C}}\times_{P,\tau} F$ is Fano if and only if  $F$ is Fano and the condition (1.1) in \cite{PS}  holds.

 Now let $X=G^{\mathbb{C}}\times_{P,\tau} F$ be Fano. As in \cite{PS}, we denote the set of $G \times T^m$-invariant (resp. $T^m$-invariant) 2-forms in $c_1(X)$ (resp. $c_1(F)$) by $c_1(X)^{G \times T^m}$ (resp. $c_1(F)^{T^m}$).  By Lemma 5.1 in \cite{PS}, the  map $R: c_1(X)^{G \times T^m}\rightarrow  c_1(F)^{T^m}$ given by $R(\tilde{\omega})=\tilde{\omega}|_{TF}$ (restriction) is invertible. We'll denote $R(\tilde{\omega})$ (not to be confused with the scalar curvature) by $\omega$. Let $E: c_1(F)^{T^m} \rightarrow  c_1(X)^{G \times T^m}$ be the inverse of the map $R$. For $\omega\in c_1(F)^{T^m}$, $E(\omega)$ is an extension of $\omega$; we'll denote $E(\omega)$ by $\tilde{\omega}$. By Lemma 5.1 in \cite{PS} $E(\omega)$ is K\"{a}hler if and only if $\omega$ is K\"{a}hler.

Fix a K\"{a}hler form $\omega_0 \in c_1(F)^{T^m}$. As in Section 2.1 there is  $u_0 \in C^\infty(F_{\text{reg}})^{T^m}$  such that $(\sqrt{-1})^m e^{-u_0}\prod_{i=1}^m(dw_i\wedge d\bar{w}_i)$ extends to a volume form on $F$ and $\omega_0|_{F_{\text{reg}}} =\frac{1}{4\pi}dd^c u_0$. For any other K\"{a}hler form $\omega \in c_1(F)^{T^m}$, we may write $\omega=\omega_0+\frac{1}{4\pi}dd^c\varphi$, where $\varphi \in C^\infty(F)^{T^m}$.
From Facts 2 and 5 in Section 2.1 we see that  the map $(\frac{1}{2\pi}\frac{\partial (u_0+\varphi)}{\partial t_1},\cdot\cdot\cdot, \frac{1}{2\pi}\frac{\partial (u_0+\varphi)}{\partial t_m})$ is the restriction to $F_{\text{reg}}$ of the metrically normalized moment map relative to  $\omega$ (under  the basis $\{d\tau_e(Z_i)^*\}$ of $\mathfrak{t}^*$), and its image ($\Delta_F  \setminus \partial \Delta_F$) does not depend on the choice of $\omega$. (Note that our coordinate system is slightly different from the one used in \cite{PS}.)
As in \cite{PS} we denote by $\mathcal{A}^{-1}$ the product
 \begin{equation*}
\prod_{\alpha \in R_\mathfrak{m}^+}(\frac{a_\alpha^i}{2\pi}\frac{\partial (u_0+\varphi)}{\partial t_i}+b_\alpha),
\end{equation*}
where  $a_\alpha^i:=\alpha(\sqrt{-1}Z_i)$ and $b_\alpha:=\alpha(\sqrt{-1}Z_V)$. Then one sees that
\begin{equation*}
\mathcal{A}^{-1}=\prod_{\alpha \in R_\mathfrak{m}^+} \sqrt{-1}\alpha(\sum_{i=1}^mf_iZ_i+Z_V),
\end{equation*}
where $f_i=\frac{1}{2\pi}\frac{\partial (u_0+\varphi)}{\partial t_i}$.  By the Fano condition on $X$, the property of the image of $(f_1, \cdot\cdot\cdot, f_m)$ noticed above, and  Theorem A in \cite{PS},  it follows  that $\mathcal{A}$ is bounded above and below by two positive constants not depending on the choice of $\omega$,
\begin{equation}\label{1}
K_1 \leq \mathcal{A} \leq K_2,
\end{equation}
compare the proof of Lemma 5.3 in \cite{PS}.

\vspace *{0.4cm}

\noindent {\bf 2.3 K\"{a}hler-Ricci flow and its reduction to PMAE}

\vspace*{0.4cm}

 Suppose  $X$ is a Fano manifold of complex dimension $n$,  and let  $g_0$   be a K\"{a}hler metric on $X$ such that  the corresponding K\"{a}hler form $\omega_0$ represents the first Chern class $c_1(X)$.
Let $h_0$ be a Ricci potential of $\omega_0$, so that
\begin{equation*}
\text{Ric}(\omega_0)-\omega_0=\frac{\sqrt{-1}}{2\pi}\partial \bar{\partial}h_0=\frac{1}{4\pi}dd^ch_0.
\end{equation*}
$h_0$ is uniquely determined up to an additive constant.

Consider the normalized K\"{a}hler-Ricci flow (NKRF) on $X$,
\begin{equation}\label{2}
\frac{\partial \omega}{\partial t}=-\text{Ric}(\omega)+\omega,    \hspace{4mm} \omega(0)=\omega_0.
\end{equation}

The following result is well-known, cf. for example \cite{TZ3}.

\begin{prop} \label{prop 2.1} \ \  Let $X$ be Fano. Fix a K\"{a}hler form $\omega_0\in c_1(X)$. A smooth family  of K\"{a}hler forms $\omega(t)$ ($t\in [0,T)$) on $X$ solves the NKRF (\ref{2}) if and only if there is a smooth family of smooth functions $\varphi(t) ( t\in [0,T))$ on $X$
with $\omega(t)=\omega_0+\frac{\sqrt{-1}}{2\pi}\partial \bar{\partial}\varphi (t)$ such that
\begin{equation}\label{3}
\frac{\partial}{\partial t}\varphi=\log \frac{\det(g_{i\bar{j}}^0+\varphi_{i\bar{j}})}{\det (g_{i\bar{j}}^0)}+\varphi-h_0,    \hspace{4mm} \varphi(0)=0, \hspace{4mm}  \omega_0+\frac{\sqrt{-1}}{2\pi}\partial \bar{\partial}\varphi (t)>0,
\end{equation}
where $(g_{i\bar{j}}^0)$ is the K\"{a}hler metric corresponding to $\omega_0$, and $h_0$ is a Ricci potential of $\omega_0$.
\end{prop}
\noindent {\bf Proof}.  It is straightforward to see that if $\varphi(t)$ is a solution of (\ref{3}), then $\omega(t):=\omega_0+\frac{\sqrt{-1}}{2\pi}\partial \bar{\partial}\varphi (t)$  solves (\ref{2}).

Conversely, suppose $\omega(t)$ solves (\ref{2}). Write $\omega(t)=\omega_0+\frac{\sqrt{-1}}{2\pi}\partial \bar{\partial}\tilde{\varphi}(t)$ using $\partial \bar{\partial}$-lemma, where $\tilde{\varphi}(t)$ is a smooth family of smooth functions on $X$ with $\tilde{\varphi}(0)=0$.  Then
\begin{equation*}
\frac{\sqrt{-1}}{2\pi}\partial \bar{\partial}\frac{\partial \tilde{\varphi}}{\partial t}=\frac{\sqrt{-1}}{2\pi}\partial \bar{\partial}(\log \frac{\det(g_{i\bar{j}}^0+\tilde{\varphi}_{i\bar{j}})}{\det (g_{i\bar{j}}^0)}+\tilde{\varphi}-h_0).
\end{equation*}
By the maximum principle there exist constants $C(t)$ (smoothly depending on $t$) such that
\begin{equation*}
\frac{\partial \tilde{\varphi}(t)}{\partial t}=\log \frac{\det(g_{i\bar{j}}^0+\tilde{\varphi}_{i\bar{j}}(t))}{\det (g_{i\bar{j}}^0)}+\tilde{\varphi}(t)-h_0+C(t).
\end{equation*}
Choose $b(t)$ such that
\begin{equation*}
\frac{d }{dt}b(t)=b(t)-C(t), \hspace{4mm} b(0)=0.
\end{equation*}
Let $\varphi(t,\cdot)=\tilde{\varphi}(t,\cdot)+b(t)$.  Then $\varphi(t)$ solves (\ref{3}).   \hfill{$\Box$}

\vspace*{0.4cm}

For $(X,\omega_0)$ as above, Cao \cite{C} showed that the solution of (\ref{2}) (or equivalently (\ref{3})) exists for $t\in [0,\infty)$.

Let $\varphi_t$ be a solution of (\ref{3}), and $\omega_t=\omega(t)=\omega_0+\frac{\sqrt{-1}}{2\pi}\partial \bar{\partial}\varphi (t)$. For each $t$ we choose a constant $c_t$ such that
\begin{equation}\label{a}
\int_Xe^{-\frac{\partial \varphi}{\partial t}+c_t}\omega_t^n=\int_X\omega_0^n.
\end{equation}
By a deep estimate of Perelman (see \cite{ST}), there is a constant $C$ depending only on $\omega_0$ such that for all $t\in [0,\infty)$,
\begin{equation}\label{4}
|-\frac{\partial \varphi}{\partial t}+c_t|\leq C.
\end{equation}

\vspace *{0.4cm}

Now we assume that the Fano manifold $X$ is a homogeneous toric bundle of the form $G^{\mathbb{C}}\times_{P,\tau} F$ as above. Consider the normalized K\"{a}hler-Ricci flow  on $X$,
\begin{equation}\label{5}
\frac{\partial \tilde{\omega}}{\partial t}=-\text{Ric}(\tilde{\omega})+\tilde{\omega},    \hspace{4mm} \tilde{\omega}(0)=\tilde{\omega}_0.
\end{equation}

\begin{prop} \label{prop 2.2} \ \  Let  $X=G^{\mathbb{C}}\times_{P,\tau} F$ be Fano. Fix a K\"{a}hler form $\tilde{\omega}_0\in c_1(X)^{G \times T^m}$ and write  $R(\tilde{\omega}_0)=\omega_0\in c_1(F)^{T^m}$. A smooth family  of K\"{a}hler forms $\tilde{\omega}(t)$ ($t\in [0,T)$) on $X$ solves the NKRF (\ref{5}) if and only if there is a smooth family of functions $\varphi(t)\in C^\infty(F)^{T^m} ( t\in [0,T))$
with $\varphi(0)=0$ and  $\omega_0+\frac{\sqrt{-1}}{2\pi}\partial \bar{\partial}\varphi (t)>0$  such that $\tilde{\omega}(t)=E(\omega_0+\frac{\sqrt{-1}}{2\pi}\partial \bar{\partial}\varphi (t))$,  and at all points of $F_{\text{reg}}$,
\begin{equation}\label{6}
\frac{\partial (u_0+\varphi)}{\partial t}=\log[\det(\frac{\partial^2 (u_0+\varphi)}{\partial t_i \partial t_j})\prod_{\alpha \in R_\mathfrak{m}^+}(\frac{a_\alpha^i}{2\pi}\frac{\partial (u_0+\varphi)}{\partial t_i}+b_\alpha)]+u_0+\varphi,
\end{equation}
where $u_0 \in C^\infty(F_{\text{reg}})^{T^m}$   such that $(\sqrt{-1})^m e^{-u_0}\prod_{i=1}^m(dw_i\wedge d\bar{w}_i)$ extends to a volume form on $F$ and $\omega_0|_{F_{\text{reg}}}=\frac{\sqrt{-1}}{2\pi}\partial \bar{\partial}u_0$.
\end{prop}

\noindent {\bf Proof}.   Let $\varphi(t)\in C^\infty(F)^{T^m} ( t\in [0,T))$ be a smooth family of functions  with $\varphi(0)=0$ and $\omega_0+\frac{\sqrt{-1}}{2\pi}\partial \bar{\partial}\varphi (t)>0$ such that  at all points of $F_{\text{reg}}$ equation (\ref{6}) is satisfied. Write  $\omega(t)=\omega_0+\frac{\sqrt{-1}}{2\pi}\partial \bar{\partial}\varphi (t)$. By Lemma 5.1 in \cite{PS},  $\tilde{\omega}_0:=E(\omega_0)\in c_1(X)^{G \times T^m}$,  and $\tilde{\omega}(t):=E(\omega(t))\in c_1(X)^{G \times T^m}$.
Observe that the cohomology class
\begin{equation*}
[\frac{\partial \tilde{\omega}(t)}{\partial t}]=\frac{\partial [\tilde{\omega}(t)]}{\partial t}=0.
\end{equation*}
Denote the RHS of (\ref{6}) by $-\Psi$. Take $\frac{1}{4\pi}dd^c$ of both sides of (\ref{6});
both sides of the new equation thus obtained can be naturally and uniquely extended to be smooth 2-forms on $F$. By Lemma 2.2 and Proposition 2.3 of \cite{PS}, $\frac{1}{4\pi}dd^c \Psi$ coincides with the restriction to $F_{\text{reg}}$ of
$\text{Ric}(\tilde{\omega})-\tilde{\omega}$.  Note that
\begin{equation*}
\frac{\partial \tilde{\omega}(t)}{\partial t}|_{TF}=\frac{\partial \omega(t)}{\partial t}.
\end{equation*}
Then we get that
\begin{equation*}
\frac{\partial \tilde{\omega}(t)}{\partial t}|_{TF}=(\tilde{\omega}(t)-\text{Ric}(\tilde{\omega}(t)))|_{TF},
\end{equation*}
or
\begin{equation*}
R(\frac{\partial \tilde{\omega}(t)}{\partial t}+\text{Ric}(\tilde{\omega}(t)))=R(\tilde{\omega}(t)).
\end{equation*}
 Using  Lemma 5.1 in \cite{PS} again, we see that
 \begin{equation*}
\frac{\partial \tilde{\omega}(t)}{\partial t}+\text{Ric}(\tilde{\omega}(t))=\tilde{\omega}(t),
\end{equation*}
 thus   $\tilde{\omega}(t)$ solves equation (\ref{5}).

Conversely assume that  $\tilde{\omega}(t)$ solves equation (\ref{5}). By Lemma 5.1 in \cite{PS} we have $\omega_0:=R(\tilde{\omega}_0)\in c_1(F)^{T^m}$ and $\omega(t):=R(\tilde{\omega}(t))\in c_1(F)^{T^m}$. We can write $\omega_t=\omega(t)=\omega_0+\frac{\sqrt{-1}}{2\pi}\partial \bar{\partial}\psi(t)$, where
$\psi_t=\psi(t)$ is a smooth family of functions in $C^\infty(F)^{T^m} ( t\in [0,T))$.  Let
\begin{equation*}
\Psi_t:=-\log[\det(\frac{\partial^2 (u_0+\psi_t)}{\partial t_i \partial t_j})\prod_{\alpha \in R_\mathfrak{m}^+}(\frac{a_\alpha^i}{2\pi}\frac{\partial (u_0+\psi_t)}{\partial t_i}+b_\alpha)]-(u_0+\psi_t).
\end{equation*}
Restricting equation (\ref{5}) to $F_{\text{reg}}$ and using Lemma 2.2 and Proposition 2.3 in \cite{PS} again, we see that
 \begin{equation*}
 \frac{1}{4\pi}dd^c\frac{\partial\psi}{\partial t}=\frac{1}{4\pi}dd^c(-\Psi_t).
 \end{equation*}
 It follows that $\frac{\partial\psi}{\partial t}=-\Psi_t+C(t)$.  Choose $b(t)$ as in the proof of Proposition 2.1 and let $\varphi(t,\cdot)=\psi(t,\cdot)+b(t)$, then $\varphi(t)$ solves equation (\ref{6}).  \hfill{$\Box$}

\vspace*{0.4cm}

Let  $X=G^{\mathbb{C}}\times_{P,\tau} F$ be Fano. Given an initial  K\"{a}hler form   $\tilde{\omega}_0\in c_1(X)^{G \times T^m}$ on $X$,   the function $h_0$   in (\ref{3}) and the function $u_0$ in (\ref{6}) are unique up to  additive constants. By Fact 3 in Section 2.1 the function $u_0-v^0$ can be extended to a $T^m$-invariant smooth function on $F$, so from Fact 4 in Section 2.1 we have that
\begin{equation*}
\sup_{F_{\text{reg}}}|\log\det(\frac{\partial^2 u_0}{\partial t_i \partial t_j})+u_0| < \infty.
\end{equation*}
Compare p. 327 in \cite{Z}. On the other hand, using Lemma 2.2 and Proposition 2.3 in \cite{PS},  we have
\begin{equation*}
\begin{split}
& \frac{\sqrt{-1}}{2\pi}\partial \bar{\partial}\{\log[\det(\frac{\partial^2 u_0}{\partial t_i \partial t_j})\prod_{\alpha \in R_\mathfrak{m}^+}(\frac{a_\alpha^i}{2\pi}\frac{\partial u_0}{\partial t_i}+b_\alpha)]+u_0+h_0\}            \\
=&  -(\text{Ric}(\tilde{\omega}_0)-\tilde{\omega}_0)|_{F_{\text{reg}}}+\frac{\sqrt{-1}}{2\pi}(\partial \bar{\partial}h_0)|_{F_{\text{reg}}}       \\
=& 0.
\end{split}
\end{equation*}
Recall also (\ref{1}).
\noindent It follows that
\begin{equation*}
\log[\det(\frac{\partial^2 u_0}{\partial t_i \partial t_j})\prod_{\alpha \in R_\mathfrak{m}^+}(\frac{a_\alpha^i}{2\pi}\frac{\partial u_0}{\partial t_i}+b_\alpha)]+u_0+h_0\
\end{equation*}
is a constant. So we can and will choose $h_0$ and  $u_0$ such that
\begin{equation}
\det(\frac{\partial^2 u_0}{\partial t_i \partial t_j})\prod_{\alpha \in R_\mathfrak{m}^+}(\frac{a_\alpha^i}{2\pi}\frac{\partial u_0}{\partial t_i}+b_\alpha)=e^{-h_0-u_0}.
\end{equation}

 From Propositions 2.1 and 2.2  we have a correspondence between the solutions to the equation (2.3) and the solutions to the equation (2.7), where $h_0$ and $u_0$ satisfy the constraint (2.8). More precisely, for a homogeneous toric bundle $X$ which is Fano with an initial K\"{a}hler form in $c_1(X)^{G \times T^m}$,  given a solution  $\varphi$ to equation (\ref{3}), $u_0+\varphi|_{F_{\text{reg}}}$ solves equation (\ref{6}), where we assume that $u_0$ satisfies the additional constraint (2.8); conversely,  given a  smooth family of functions $\varphi(t)\in C^\infty(F)^{T^m}$ with $\varphi(0)=0$ and  $\omega_0+\frac{\sqrt{-1}}{2\pi}\partial \bar{\partial}\varphi (t)>0$   such that at all points of $F_{\text{reg}}$, $u_0+\varphi$ solves  equation (\ref{6}),  the $G\times T^m$-invariant extension of $\varphi$ to the whole $X$ (which will also be denoted by $\varphi$) solves equation (\ref{3}), where $h_0$ satisfies the additional constraint (2.8).

The equation (\ref{6}) with $\varphi(0)=0$ is actually a parabolic real Monge-Amp\`{e}re equation,
\begin{equation}\label{7}
\frac{\partial (u_0+\varphi)}{\partial t}=\log\det(\frac{\partial^2 (u_0+\varphi)}{\partial t_i \partial t_j})+u_0+\varphi-\log \mathcal{A}, \hspace{4mm}   on \hspace{2mm} \mathbb{R}^m,
\end{equation}
with $\varphi(0)=0$ (and the matrices
$(\frac{\partial^2 (u_0+\varphi)}{\partial t_i \partial t_j})>0$), where $ \mathcal{A}$ is defined in Section 2.2.

\vspace*{0.4cm}

\section{ Proof of Theorem 1.1}

We will write $\underline{t}=(t_1,\cdot\cdot\cdot, t_m)\in \mathbb{R}^m$, and $d\underline{t}=dt_1\cdot\cdot\cdot dt_m$.

\begin{lem} \label{3.1} (\cite{WZ}, see also  Proposition 2 on p. 54 in \cite{Do})  Assume that $v$ is a smooth convex function on $\mathbb{R}^m$ which attains minimal value 0, such that $\det (D^2v)\geq \lambda$ when $v\leq 1$. Let $\Gamma$ be the set where $v\leq 1$. Then $\text{Vol}(\Gamma) \leq C\lambda^{-1/2}$, where $C$ is a constant depending only on the dimension $m$.
\end{lem}
\noindent {\bf Proof}.   We follow the outline  given on p. 54 of \cite{Do}.
 $\Gamma$ is a bounded convex set in $\mathbb{R}^m$ with nonempty interior. Let  $E$ be the ellipsoid of minimum volume containing $\Gamma$ centered at the barycenter of $\Gamma$. By a variant of John's Lemma (see for example Theorem 1.8.2 in \cite{Gu}), $\alpha_mE\subset \Gamma \subset E$, where $\alpha_m=m^{-3/2}$ and $\beta Y$ denotes the $\beta$-dilation of a bounded set $Y\subset \mathbb{R}^m$ w.r.t. its barycenter, that is $\beta Y:=\{\underline{t_0}+\beta(\underline{t}-\underline{t_0}) \hspace{1mm} | \hspace{1mm} \underline{t}\in Y\}$, where $\underline{t_0}$ is the barycenter of $Y$.
Choose a unimodular affine transformation $T$ such that $B(0,\alpha_mr) \subset T(\Gamma) \subset B(0,r)$.  We write $T\underline{t}=A\underline{t}+b$ for any $\underline{t}\in \mathbb{R}^m$, where $A$ is a $m \times m$ matrix with $\det A=1$ and $b\in \mathbb{R}^m$.
Let
\begin{equation*}
\tilde{v}(\underline{t})=v(A^{-1}(r\underline{t}-b))-1.
\end{equation*}
Let $\Gamma'=\frac{1}{r}T(\Gamma)$.
   Then the barycenter of $\Gamma'$ is the origin 0, $B(0,\alpha_m) \subset \Gamma' \subset B(0,1)$,   $\min_{ \Gamma'}\tilde{v}=-1$ and $\tilde{v}=0$ on $\partial \Gamma'$.
We also have $\det D^2 \tilde{v}\geq r^{2m}\lambda$ on $\Gamma'$.
From the proof of Proposition 3.2.4 in \cite{Gu} we have
\begin{equation*}
\int_{\frac{1}{2}\Gamma'} \det D^2 \tilde{v} d\underline{t} \leq C |\min_{ \frac{1}{2}\Gamma'}\tilde{v}|^m\leq C,
\end{equation*}
where $C$  depends only on the dimension $m$.  It follows that
\begin{equation*}
r^{2m}\lambda \text{Vol} (\frac{1}{2}\Gamma') \leq C,
\end{equation*}
 and
\begin{equation*}
\frac{r^m\lambda}{2^m} \text{Vol}(\Gamma) = \frac{r^m\lambda}{2^m}\text{Vol} (T(\Gamma))\leq C.
 \end{equation*}
  We also have
 \begin{equation*}
 \text{Vol} (\Gamma)= \text{Vol} (T(\Gamma)) \leq \omega_mr^m,
   \end{equation*}
  where $\omega_m$ is the volume of the unit $m$-ball.
It follows that
\begin{equation*}
\frac{\lambda}{2^m}\frac{\text{Vol}(\Gamma)}{\omega_m}\text{Vol}(\Gamma)\leq C,
 \end{equation*}
and we are done.         \hfill{$\Box$}

\vspace*{0.4cm}

For  a homogeneous toric bundle $X$ which is Fano with an initial K\"{a}hler form $\omega_0 \in c_1(X)^{G \times T^m}$,  given a solution  $\varphi$ to equation (\ref{3}), from above we see that $u=u_t=u(t,\cdot)=u_0+\varphi$ is a solution of equation (\ref{7}) with $\varphi(0)=0$, where $u_0$ is as in Section 2.3.

The following result is well known.

\begin{lem} \label{lem 3.2} Let $u$ be as above. Then
\begin{equation*}
\int_{\mathbb{R}^m} \det D^2 u_t d\underline{t}=\text{Vol}(\Delta).
\end{equation*}
\end{lem}
\noindent {\bf Proof}. This follows immediately from Fact 2 and the change of variable formula for multiple integrals.
\hfill{$\Box$}

\begin{lem} \label{lem 3.3}  Let $u$ be as above and $\bar{u}=\bar{u}_t=u_t-c_t$, where $c_t$ is the  constant determined by equation (\ref{a}) with $\omega_t$ and  $ \omega_0$ there replaced by $\tilde{\omega}_t$ and $ \tilde{\omega}_0$ respectively. Let $m_t=\inf_{\mathbb{R}^m}\bar{u}_t(\underline{t})$. Then there is a constant $C$ such that for any $t\in (0,\infty)$,
\begin{equation*}
|m_t|\leq C.
\end{equation*}
\end{lem}
\noindent {\bf Proof}.  We follow \cite{Z} to argue that $m_t$ is bounded below.
As noted above,  $\frac{1}{2\pi}Du=(\frac{1}{2\pi}\frac{\partial (u_0+\varphi)}{\partial t_1}, \cdot\cdot\cdot, \frac{1}{2\pi}\frac{\partial (u_0+\varphi)}{\partial t_m}): \mathbb{R}^m \rightarrow \Delta_F \setminus \partial \Delta_F$. So
\begin{equation*}
|D\bar{u}_t|=|Du_t|\leq 2\pi \text{diam}(\Delta_F).
\end{equation*}
Also note  that  since $0\in  \Delta_F \setminus \partial \Delta_F$ and $\bar{u}_t$ is convex, $\bar{u}_t$ attains $m_t$.

 From (\ref{7}) we have
\begin{equation*}
\int_{\mathbb{R}^m}e^{-\bar{u}}d\underline{t}=\int_{\mathbb{R}^m}\mathcal{A}^{-1}e^{-\frac{\partial u}{\partial t}+c_t} \det D^2u d\underline{t},
\end{equation*}
combining this with  (\ref{1}), (\ref{4}) and Lemma 3.2  one sees that there are positive constants $c_1$ and $c_2$ independent of $t$, such that
\begin{equation*}
c_1 \leq \int_{\mathbb{R}^m}e^{-\bar{u}}d\underline{t} \leq c_2.
\end{equation*}
Then  one easily sees  that $m_t$ is uniformly bounded below.

To show $m_t$ is bounded above we follow the argument in Section 3.4 of \cite{Do}, which is a slight variant of that in \cite{WZ}.
Let $v=v_t=\bar{u}_t-m_t$. From (\ref{1}), (\ref{4}) and (\ref{7}) we get
\begin{equation*}
  \det D^2v
= \det D^2 u
= \exp \{\frac{\partial u}{\partial t}-c_t-\bar{u}+\log \mathcal{A}\}
 \geq  ce^{-\bar{u}},
\end{equation*}
where $c$ is a positive constant independent of $t$. It follows
\begin{equation*}
\det D^2v_t \geq \frac{c}{e}e^{-m_t} \hspace{4mm}  \text{when}  \hspace{4mm}  v_t\leq 1.
\end{equation*}
Let $\Gamma$ be the set where $v_t\leq 1$. By Lemma 3.1,
\begin{equation*}
\text{Vol}(\Gamma)\leq C (\frac{c}{e}e^{-m_t})^{-1/2}=C'e^{m_t/2}.
\end{equation*}
For each $s>0$ let $\Gamma_s$ be the set $\{v_t \leq s\}$ and $V(s)=\text{Vol}(\Gamma_s)$. By convexity $\Gamma_s$ is contained in the $s$-dilation of $\Gamma$ w.r.t. the minimum point of $v_t$.  Then
\begin{equation}
V(s)\leq s^m\text{Vol}(\Gamma) \leq C's^me^{m_t/2}.
\end{equation}
Note that from (\ref{1}), (\ref{4}) and (\ref{7}) we have
\begin{equation}
\det D^2u_t\leq C e^{-m_t}e^{-v_t}.
\end{equation}
  So by Lemma 3.2 and (3.2) we have
\begin{equation}
C_1 = \int_{\mathbb{R}^m} \det D^2 u_t d\underline{t} \leq Ce^{-m_t} \int_{\mathbb{R}^m}e^{-v_t}d\underline{t}.
\end{equation}

Using the co-area formula
\begin{equation*}
\int_{\mathbb{R}^m}e^{-v}d\underline{t}=\int_0^\infty e^{-s}V(s)ds
\end{equation*}
and (3.1), we get from (3.3) that
\begin{equation*}
C_1 \leq  Ce^{-m_t}\int_0^\infty     C' e^{-s}s^me^{m_t/2} ds =  C''e^{-m_t/2},
\end{equation*}
and we are done.
\hfill{$\Box$}

\vspace*{0.4cm}

  Let $x_t$ be the minimal point of $\bar{u}_t$, and $\bar{\bar{u}} =\bar{\bar{u}}_t(\cdot)=\bar{u}_t(\cdot+x_t)-m_t$.  Let $\bar{\varphi}(\cdot)=\bar{\bar{u}}(\cdot)-u_0(\cdot)$.  We extend $\bar{\varphi}$ to  a  $G \times T^m$-invariant function on $X$, still  denoted by  $\bar{\varphi}$. Using Facts 1, 2 and 3 in Section 2.1, convexity of $\bar{\bar{u}}$, (2.1), Perelman's estimate (2.5), and Lemma 3.3, one can show that
\begin{equation}
|\sup \bar{\varphi}| \leq C.
\end{equation}
Compare the proof of Lemma 3.4 in \cite{WZ} and Proposition 3.2 in \cite{Z}.

Using Proposition 2.1 in \cite{TZ2}, Lemma 4.4 in \cite{PS}, the monotonicity (see (4.5) in \cite{TZ3}) of the generalized K-energy $\tilde{\mu}$ (introduced in \cite{TZ2}) along certain modified K\"{a}hler-Ricci flow, Lemma 5.1 in \cite{TZ2}, (3.4) above, Lemma 3.1 in \cite{CTZ}, Proposition 3.1 in \cite{TZ3},  and Lemma 3.3,  one can show that
\begin{equation}
|| \bar{\varphi}||_{C^0(X)} \leq C.
\end{equation}
Compare the proof of Proposition 4.1 in \cite{Z}.  Using Lemma 5.1 in \cite{TZ2}, Lemma 3.3 and (3.5) one can show that
\begin{equation}
\tilde{\mu}(\varphi) \geq -C,
\end{equation}
where $\varphi$ is a solution of (2.3) as above; compare the proof of Corollary 4.2 in \cite{Z}.
Then by  using the monotonicity  of the generalized K-energy   along certain modified K\"{a}hler-Ricci flow  and (3.6), one can show that
\begin{equation*}
|\frac{\partial \varphi}{\partial t}| \leq C,
\end{equation*}
after suitably normalizing the Ricci potential $h_0$ of the initial metric.  Compare Proposition 4.3 in \cite{Z}.

 As in Section 5 of \cite{Z} we can find $x_t' \in \mathbb{R}^m$ such that
 \begin{equation*}
|x_t-x_t'|\leq C   \hspace{0.2cm} \text{and}    \hspace{0.2cm} |\frac{dx_t'}{dt}| \leq C.
\end{equation*}
 Let $\tilde{u}(\cdot)=\tilde{u}_t(\cdot)=u_t(\cdot + x_t')$ and $\tilde{\varphi}_t(\cdot)=\tilde{u}_t(\cdot)-u_0(\cdot)$.  The $G \times T^m$-invariant extension of $\tilde{\varphi}_t$ gives a K\"{a}hler potential on $X$ relative to $\omega_0$, still denoted by $\tilde{\varphi}_t$. Moreover, $\frac{dx_t'}{dt}$ are corresponding to a family of holomorphic vector fields on $X$, denoted by $\tilde{Y}_t$. (The real part of) $\tilde{Y}_t$ generate a family of elements, denoted by $\rho_t$, in the algebraic torus subgroup $(T^m)^\mathbb{C}$ of Aut($X,J$). Then $\rho_t^*\omega_{\varphi_t}=\omega_{\tilde{\varphi}_t}$, where  $\omega_{\varphi_t}=\omega_0+\frac{\sqrt{-1}}{2\pi}\partial \bar{\partial}\varphi_t$ and  $\omega_{\tilde{\varphi}_t}=\omega_0+\frac{\sqrt{-1}}{2\pi}\partial \bar{\partial}\tilde{\varphi}_t$.  As in \cite{Z}, using the  estimates above, we can show that
   \begin{equation*}
||\tilde{\varphi}_t||_{C^0(X)}\leq C   \hspace{0.2cm} \text{and}    \hspace{0.2cm} |\frac{\partial \tilde{\varphi}}{\partial t}| \leq C.
\end{equation*}

 With the above preparation we can proceed  as in Section 5 of \cite{Z}, using arguments as in Section 6 of \cite{TZ3}  and the uniqueness theorem in \cite{TZ1} \cite{TZ2}, and get that the K\"{a}hler metrics $\omega_{\tilde{\varphi}_t}$ converge to a K\"{a}hler-Ricci soliton as $t\rightarrow \infty$.

\vspace*{0.4cm}

\noindent {\bf Acknowledgements}.   I'm grateful to Professor Xiaohua Zhu for answering my questions on his paper [30]. I would also like to thank Professor Bin Zhou  for helpful discussions.
I also thank the referee for the comments. I'm partially supported by Beijing Natural Science Foundation (Z190003).


\hspace *{0.8cm}

\vspace *{0.4cm}

School of Mathematical Sciences, Beijing Normal University,

Laboratory of Mathematics and Complex Systems, Ministry of Education,

Beijing 100875, P.R. China

 E-mail address: hhuang@bnu.edu.cn

\end{document}